\documentclass[11pt,a4paper]{article}
\usepackage{amsfonts}
\usepackage{}

\usepackage{graphicx}
\usepackage{amsfonts,amsmath, amssymb}

\newtheorem{theo}{\textbf{\ \ \quad Theorem}}[section]

\newtheorem{remark}{\textbf{\ \ \quad Remark}}[section]

\newtheorem{prop}{\textbf{\ \ \quad Proposition}}[section]

\newcommand{\lbl}[1]{\label{#1}}

\newcommand{\be}{\begin{equation}}
\newcommand{\ee}{\end{equation}}
\newcommand\bes{\begin{eqnarray}}
\newcommand\ees{\end{eqnarray}}
\newcommand{\bess}{\begin{eqnarray*}}
\newcommand{\eess}{\end{eqnarray*}}

\newcommand{\nm}{\nonumber}

{\setlength\arraycolsep{2pt}}

\title{Impact of noise on Parabolic equations}

\author{Guangying Lv$^{a,b}$,  Jinlong Wei$^c$\\
\\
\ \\
   {\small \it $^a$ Institute of Applied Mathematics, Henan University, Kaifeng, Henan 475001, China}\\
   {\small \it $^c$ School of Statistics and Mathematics, Zhongnan University of}\\
   {\small \it Economics and Law, Wuhan 430073, China}\\
    {\small \tt  weijinlong.hust@gmail.com }}

\begin{document}
\maketitle

\medskip

\begin{abstract}
In this short paper, we focus on the blowup phenomenon of stochastic parabolic equations.
We first consider the probability of the event that the solutions keep positive.
Then, the blowup phenomenon in the whole space is considered. The probability of
the event that the solutions blow up in finite time is given. Lastly,
we obtain the probability of the event that blowup time of stochastic parabolic
equations large than or less than the deterministic case.

{\bf Keywords}: Blowup; Stochastic heat equation; Impact of noise.

AMS subject classifications (2010): 35K20, 60H15, 60H40.

\end{abstract}

\baselineskip=15pt

\section{Introduction}
\setcounter{equation}{0}

For a deterministic partial differential equations, when we add a noise on it,
we first want to know how to change about the solution, that is, the effect of noise.
More precisely, if
the solutions of deterministic parabolic equations keep positive, what
probability of the solutions keep positive for the stochastic case ? In this paper, we will give part of
positive answer. Similarly, for the blowup phenomenon, we want to know
the probability of the event that the solutions blow up in finite time.

We firstly recall some known results of stochastic partial differential equations (SPDEs).
In this paper, we only focus on the stochastic parabolic equations.
It is known that the existence and uniqueness of
global solutions to SPDEs can  be established under appropriate conditions (\cite{Cb2007}).
For the finite time blowup phenomenon of stochastic parabolic equations, we first consider
the case on a bounded domain. Consider the following equation
   \bes\left\{\begin{array}{lll}
   du=(\Delta u+f(u))dt+\sigma(u)dW_t, \ \ \qquad t>0,&x\in  D,\\[1.5mm]
   u(x,0)=u_0(x)\geq0, \ \ \ &x\in D,\\[1.5mm]
   u(x,t)=0, \qquad \qquad \qquad \qquad \qquad \qquad t>0,  &x\in\partial D,
    \end{array}\right.\lbl{1.2}\ees
Da Prato-Zabczyk \cite{PZ1992} studied  the existence of global solutions of (\ref{1.2})  with
constant $\sigma$. Manthey-Zausinger \cite{MZ1999} considered (\ref{1.2}),
with $\sigma$ satisfying a global Lipschitz condition. Dozzi and L\'{o}pez-Mimbela
\cite{DL2010} studied  equation (\ref{1.2}) with
$\sigma(u)=u$ and  proved that if $f(u)\geq u^{1+\alpha}$ ($\alpha>0$) and the initial data is large enough, the solution will blow up in finite time,
and that if $f(u)\leq u^{1+\beta}$ ($\beta$ is a certain positive constant) and the initial data is
small enough, the solution will exist globally (also see \cite{NX2012}). When
$\sigma$ does not satisfy the global Lipschitz condition,
Chow \cite{C2009,C2011} obtained the finite time blowup phenomenon.
Lv-Duan \cite{LD2015} described the competition between the nonlinear term and noise term for equation (\ref{1.2}).
Bao-Yuan \cite{BY2014} and Li et al.\cite{LPJ2016} obtained the existence of local solutions of (\ref{1.2}) with jump process
and L$\acute{e}$vy process, respectively. For blowup phenomenon of stochastic functional
parabolic equations, see \cite{CL2012,LWW2016} for details.

We remark that the method used to prove the finite time blowup on bounded domain is
the stochastic Kaplan's first eigenvalue method. In order to make sure the inner product
$(u,\phi)$ is positive, the authors firstly proved the solutions of (\ref{1.2}) keep positive under some
assumptions, see \cite{BY2014,C2009,C2011,LPJ2016,LD2015}. The method used to prove the positivity of
solutions is that the negative part is zero. The main difficulty is to choose suitable
test function. In the present paper, we will introduce another method to prove the positivity
of solutions.
We also remark that, in our paper \cite{LW2019}
a new method (stochastic concavity method)  is introduced to prove the solutions blow up in finite time.
The advantage of this method is that we need not the positivity of solutions.

In former papers \cite{FLN2018,LW2019}, the blowup phenomenon in the whole space is considered
in the form of $\mathbb{E}u^2(x,t)$. That is to say, the moment of the solutions will
blow up in finite time. From the point of probability theory, we want to know the probability
of event that the solutions blow up in finite time. In present paper, we study the
parabolic equations with linear multiplicative noise and give the probability of the event.

On the other hand, we remark that the existence of finite time blowup solution was obtained
by Dozzi and L\'{o}pez-Mimbela
\cite{DL2010}. But the estimate of blowup time is no result. This is our
second aim. We will estimate the probability that blowup time of stochastic parabolic
equations large than or less than the deterministic case.

The advantage of linear multiplicative noise is that we can change stochastic parabolic equations into
 random parabolic equations. And then we can use the comparison principle and the results of
deterministic case to get the results of the stochastic case.

Throughout this paper, we write $C$ as a general positive constant and $C_i$, $i=1,2,\cdots$ as concrete positive constants.

\section{The impact of additive noise}
\setcounter{equation}{0}
In this section, we consider the impact of additive noise on parabolic equation.
Our aim is to find the probability of the event that the solutions keep positive or
belong to some interval or are lower (larger) than the solutions of the corresponding deterministic case.

We first consider a simple case:
    \bes \left\{\begin{array}{llll}
du_t=\Delta udt+\sigma dB_t,\ \ t>0,\ &x\in\mathbb{R}^d,\\
u(x,0)=u_0(x)\gneqq0, \ \  &x\in\mathbb{R}^d,
   \end{array}\right.\lbl{2.1}\ees
where $\sigma>0$, and $B_t$ is one-dimensional Brownian motion.
A mild solution to (\ref{2.1}) in sense of Walsh \cite{walsh1986} is any
$u$ which is adapted to the filtration generated by the
white noise and satisfies the following evolution equation
   \bess
u(x,t)=\int_{\mathbb{R}^d}K(x-y,t)u_0(y)dy
+\int_0^t\int_{\mathbb{R}^d}K(x-y,t-s)\sigma dydB_s,
   \eess
where $K(x,t)$ denotes the heat kernel of Laplacian operator, i.e.,
   \bess
K(x,t)=\frac{1}{(4\pi t)^{\frac{d}{2}}}\exp\left(-\frac{|x|^2}{4t} \right)
  \eess
satisfies
   \bess
\left(\frac{\partial}{\partial t}-\Delta\right)K(x,t)=0\ \ \ {\rm for}\ \
(x,t)\neq(0,0).
   \eess
Due to the properties of heat kernel $K$, we have
\bess
u(x,t)=\int_{\mathbb{R}^d}K(x-y,t)u_0(y)dy
+\sigma  B_t,
   \eess
which implies that
    \bess
\mathbb{P}(u(x,t)>0)=\mathbb{P}\Big(B_t>\frac{A_t(x)}{\sigma\sqrt{t}}\Big)
=1-\Phi\Big(\frac{A_t(x)}{\sigma\sqrt{t}}\Big),
   \eess
where $A_t(x)=-\int_{\mathbb{R}^d}K(t,x-y)u_0(y)dy$.
Similarly, we have $\mathbb{P}(u(x,t)\leq0)=\Phi\left(\frac{A_t(x)}{\sigma\sqrt{t}}\right)$.
Here $\Phi(x)$ is the distribution function of normal random variable.

Similarly, for $a,b\in\mathbb{R}$ and $a<b$, we have
   \bess
\mathbb{P}(a<u(x,t)\leq b)&=&\mathbb{P}\left(\frac{a+A_t(x)}{\sigma\sqrt{t}}<\frac{B_t}{\sqrt{t}}\leq\frac{b+A_t(x)}{\sigma\sqrt{t}}\right)\\
&=&\Phi\left(\frac{b+A_t(x)}{\sigma\sqrt{t}}\right)-\Phi\left(\frac{a+A_t(x)}{\sigma\sqrt{t}}\right).
   \eess
Therefore, we have the following results. Here $C-$ means that the constants is a little lower than $C$, i.e.,
$C- >C-\varepsilon$ for any $0<\varepsilon\ll1$.
   \begin{theo}\lbl{t2.1}
Assume that the initial data $u_0\geq0$ is a bounded continuous function.
Then the solution of (\ref{2.1}) will keep positive with the probability
$1-\Phi\left(\frac{A_t(x)}{\sigma\sqrt{t}}\right)$. For real numbers $a<b$,
we have
      \bess
\mathbb{P}(a<u(x,t)\leq b)=\Phi\left(\frac{b+A_t(x)}{\sigma\sqrt{t}}\right)-\Phi\left(\frac{a+A_t(x)}{\sigma\sqrt{t}}\right).
   \eess
Moreover, letting $\sigma\to0$, the probability converges to $1$ exponentially, and we get the exact rate. That is to say,
the probability
   \bess
\mathbb{P}(u(x,t)>0)\to1,\ \ {\rm as}\ \sigma\to0
    \eess
with the exponential rate $\beta$, where $\beta>0$ is any fixed constant. More precisely, we have the following
estimate
  \bess
1-\mathbb{P}(u(x,t)>0)=O(e^{-\frac{A^2_t(x)}{2\sigma^2 t}}).
  \eess
When $a+A_t(x)$ and $b+A_t(x)$ have the same sign, the event $\{a<u(x,t)\leq b\}$
will become impossible event as $\sigma\to0$.
  \end{theo}

{\bf Proof.} From the above discussion, we only note that
  \bess
1-\mathbb{P}(u(x,t)>0)=\Phi\left(\frac{A_t(x)}{\sigma\sqrt{t}}\right)
=\frac{1}{\sqrt{2\pi}}
\int_{-\infty}^{\frac{A_t(x)}{\sigma\sqrt{t}}}e^{-\frac{y^2}{2}}dy,
   \eess
and for any fixed positive constant $\delta$
   \bess
e^{-\frac{A^2_t(x)}{2\sigma^2 t}-}\int_{-\infty}^{\frac{A_t(x)}{\sigma\sqrt{t}}}e^{-\frac{y^2}{2}}dy\to0.
   \eess
The proof is complete. $\Box$

\begin{remark}\lbl{r2.1}
1. The assumption that the initial data $u_0\geq0$ can be deleted, that is to say,
for additive noise, the event that solutions are positive or negative is possible
event. In other words, the event that solutions always keep positive is not certain event.
Meanwhile, we note that if $u_0\geq0$, then $\mathbb{P}(u(x,t)\geq0)\geq\frac{1}{2}$.

2. If $a<u_0\leq b$, then $a+A_t(x)<0$ and $b+A_t(x)\geq0$, then as $\sigma\to0$,
the solutions will belong to $(a,b]$. Furthermore, we have the exact convergence rate.

3. It is easy to see that Theorem \ref{t2.1} also holds if the operator $\sigma$
is replaced by $-(-\sigma)^\alpha$ with $\alpha\in(0,1)$. More generally, if the
operator $\mathcal{A}$ has heat kernel, then Theorem \ref{t2.1} will be reasonable.
For example, we can take $\mathcal{A}=\Delta+V(\cdot)\cdot\nabla$.

4. Theorem \ref{t2.1} is similar to Large Deviation Principle, but there is a big difference from the
classical theory. We give the description about the event, i.e., how to become to the certain event.
\end{remark}

Now, we compare the solutions of stochastic parabolic with the corresponding
deterministic case.
For simplicity, we consider the following problems
    \bes \left\{\begin{array}{llll}
du_t=(\Delta u+ku)dt+\sigma(x,t) dB_t,\ \ t>0,\ &x\in\mathbb{R}^d,\\
u(x,0)=u_0(x)\gneqq0, \ \  &x\in\mathbb{R}^d,
   \end{array}\right.\lbl{4.1}\ees
and
    \bes \left\{\begin{array}{llll}
\frac{\partial}{\partial t}v=\Delta v+kv,\ \ t>0,\ &x\in\mathbb{R}^d,\\
v(x,0)=u_0(x)\gneqq0, \ \  &x\in\mathbb{R}^d,
   \end{array}\right.\lbl{4.2}\ees
where $k,\sigma>0$.
   \begin{theo}\lbl{t2.2}
Assume that the initial data $u_0$ is a bounded continuous function.
Denote the event as
   \bess
\mathcal{A}_t(x)=\{\omega\in\Omega: u(x,t,\omega)\leq v(x,t)\},
   \eess
then $\mathbb{P}(\mathcal{A}_t(x))=\frac{1}{2}$. Consequently, $\mathbb{E}u(x,t)=v(x,t)$ for any
$x\in\mathbb{R}^d$ and $t>0$.
  \end{theo}

{\bf Proof.} Let $w=u-v$, then $w$ satisfies that
    \bess \left\{\begin{array}{llll}
dw_t=(\Delta w+kw)dt+\sigma(x,t) dB_t,\ \ t>0,\ &x\in\mathbb{R}^d,\\
w(x,0)=0, \ \  &x\in\mathbb{R}^d.
   \end{array}\right.\eess
Denote $\tilde w=e^{-kt}w$, then we have
    \bes \left\{\begin{array}{llll}
d\tilde w_t=\Delta \tilde w dt+e^{-kt}\sigma(x,t) dB_t,\ \ t>0,\ &x\in\mathbb{R}^d,\\
\tilde w(x,0)=0, \ \  &x\in\mathbb{R}^d.
   \end{array}\right.\lbl{4.3}\ees
Then the solution of (\ref{4.3}) can be expressed as
   \bess
\tilde w(x,t)=\int_0^t\left(\int_{\mathbb{R}^d}K(x-y,t-s)e^{-ks}\sigma(y,s)dy\right) dB_s.
   \eess
Let $f(x,s,t)=\int_{\mathbb{R}^d}K(x-y,t-s)e^{-ks}\sigma(y,s)dy$, then
$\tilde w(x,t)=\int_0^tf(x,s,t)dB_s$. For any fixed $x\in\mathbb{R}^d$ and $t>0$, we have
$\int_0^tf(x,s,t)dB_s$ is a Guass process, whose expectation is $0$. And we also remark that
   \bess
\mathbb{P}(\mathcal{A}_t(x))=\mathbb{P}(w\leq0)=\mathbb{P}(\tilde w\leq0)=\frac{1}{2}.
   \eess
The proof is complete. $\Box$

We only focus on the linear parabolic equation in Theorems \ref{t2.1} and \ref{t2.2}. Actually, we can also consider
the nonlinear parabolic equation
         \bes \left\{\begin{array}{llll}
du_t=(\Delta u+ku^p)dt+\sigma  dB_t,\ \ t>0,\ &x\in\mathbb{R}^d,\\
u(x,0)=u_0(x)\gneqq0, \ \  &x\in\mathbb{R}^d,
   \end{array}\right. \lbl{21.1}\ees
where $k\in\mathbb{R}$.
For nonlinear parabolic equation (\ref{21.1}) with additive noise, it is impossible that
the solutions keep positive almost surely. However, we can use the Jensen's inequality
to deal with some special case.  Since the proof is easy, we only give the result and omit the proof details here.

\begin{prop}\lbl{p2.1} Assume that $p$ is an even positive number and $k>0$, then it holds that
   \bess
\mathbb{P}(u(x,t)\geq0)\geq1-\Phi\left(\frac{A_t(x)}{\sigma\sqrt{t}}\right).
   \eess
Assume that $p$ is an even positive number and $k<0$, then it holds that
   \bess
\mathbb{P}(u(x,t)\leq0)\geq\Phi\left(\frac{A_t(x)}{\sigma\sqrt{t}}\right).
   \eess
For $0<p\leq1$, we have that $\mathbb{E}|u|$ is a lower
solution of the following equation
\bess \left\{\begin{array}{llll}
\frac{\partial}{\partial t}v=\Delta v+|k|v^p,\ \ t>0,\ &x\in\mathbb{R}^d,\\
v(x,0)=u_0(x)\gneqq0, \ \  &x\in\mathbb{R}^d,
   \end{array}\right. \eess
where $u$ is a solution to (\ref{21.1}). Consequently, when $0<p\leq1$, the
solution of (\ref{21.1}) exists globally almost surely.
  \end{prop}

\section{The impact of linear multiplicative noise}
\setcounter{equation}{0}
In this section, we consider the impact of linear multiplicative noise on parabolic equations.
Our aim is to get the probability of the event that the solutions keep positive or the solutions
are lower (larger) than those of corresponding deterministic case, and so on.

Firstly, we consider the multiplicative noise.
    \bes \left\{\begin{array}{llll}
du_t=\Delta udt+ f(u)dt+\sigma u dB_t,\ \ t>0,\ &x\in\mathbb{R}^d,\\
u(x,0)=u_0(x)\gneqq0, \ \  &x\in\mathbb{R}^d,
   \end{array}\right.\lbl{2.2}\ees
where $\sigma>0$, and $B_t$ is a one-dimensional Brownian motion.
By using the It\^{o} formula, it is easy to see that if we let
$v(x,t)=e^{-\sigma B_t}u(x,t)$, then $v(x,t)$ satisfies that
    \bess\left\{\begin{array}{llll}
\frac{\partial}{\partial t}v(x,t)=\Delta v(x,t)-\frac{\sigma^2}{2}v(x,t)
+ e^{-\sigma B_t}f(e^{\sigma B_t}v),\ \ t>0,\ &x\in\mathbb{R}^d,\\
v(x,0)=u_0(x)\gneqq0, \ \  &x\in\mathbb{R}^d.
  \end{array}\right.\eess
Denote $w(x,t)=e^{\frac{\sigma^2}{2}t}v$, then we have
    \bes\left\{\begin{array}{llll}
\frac{\partial}{\partial t}w(x,t)=\Delta w(x,t)
+ e^{\frac{\sigma^2}{2}t-\sigma B_t}f(e^{\sigma B_t-\frac{\sigma^2}{2}t}w),\ \ t>0,\ &x\in\mathbb{R}^d,\\
w(x,0)=u_0(x)\gneqq0, \ \  &x\in\mathbb{R}^d.
  \end{array}\right.\lbl{2.3}\ees
Therefore, using the comparison principle, we have the following result.
\begin{theo}\lbl{t3.1} Assume that the nonlinear term $f$ satisfies the local Lipschitz, 
then the Cauchy problem {\rm(\ref{2.2})} with nonnegative initial datum admits
a local solution. Moreover, the solution 
 remains positive:
$u(x,t)\geq0$, a.s. for almost every $x\in \mathbb{R}^d$ and for all $  t\in[0,T]$, where 
$T$ is the lifetime.
\end{theo}

 \begin{remark}\lbl{r3.1} 
 Comparing Theorem \ref{t3.1} with \cite[Theorem 4.1]{LW2019}, we find the proof here
is simple and the result is exact same as the deterministic case. The reason is that 
the problem (\ref{2.3}) is random case and thus we can use the comparison principle.

Moreover, the result of Theorem \ref{t3.1} is better than that of \cite[Theorem 4.1]{LW2019}. 
More precisely, in \cite[Theorem 4.1]{LW2019}, the assumption about $f$ is that 
$f(u)\geq0$ for $u\leq0$, which stronger than that in 
this paper. 
   \end{remark}

\begin{theo}\lbl{t3.2} Assume all conditions of Theorem \ref{t3.1} hold and $u$ is
the solution of (\ref{2.2}) with $f(u)=u^p$. Then for $p>1$, $\mathbb{E}u$ is a supper
solution of the following equation
\bes \left\{\begin{array}{llll}
\frac{\partial}{\partial t}v=\Delta v+v^p,\ \ t>0,\ &x\in\mathbb{R}^d,\\
v(x,0)=u_0(x)\gneqq0, \ \  &x\in\mathbb{R}^d.
   \end{array}\right. \lbl{3.0}\ees
Consequently, when $p>1$, $\mathbb{E}u$ will blow up in finite time if
the initial data belongs to $\mathcal{U_1}$, see the following for the definition of $\mathcal{U_1}$, where
 \bess
\mathcal{U_1}=\left\{v(x)|v(x)\in BC(\mathbb{R}^d,\mathbb{R}_+),
v(x)\geq c e^{-k|x|^2},\ k>0,c\gg1\right\}.
  \eess
Here $BC=\{$bounded and uniformly continuous functions$\}$, see
Fujita \cite{F1966,F1970} and Hayakawa \cite{kH1973}.
\end{theo}

{\bf Proof.} It follows from Theorem \ref{t3.1} that $u\geq0$ almost surely.
The mild solution of (\ref{2.2}) can be expressed as
 \bess
u(x,t)&=&\int_{\mathbb{R}^d}K(t,x-y)u_0(y)dy
+\int_0^t\int_{\mathbb{R}^d}K(t-s,x-y)u^p(y,s)dyds\\
&&
+\int_0^t\int_{\mathbb{R}^d}K(t-s,x-y)u(y,s)dydB_s.
   \eess
Taking expectation in the above equality, we have
 \bess
\mathbb{E}u(x,t)\geq\int_{\mathbb{R}^d}K(t,x-y)u_0(y)dy
+\int_0^t\int_{\mathbb{R}^d}K(t-s,x-y)(\mathbb{E}u)^p(y,s)dyds,
   \eess
which implies that $\mathbb{E}u$ is a supper
solution of equation (\ref{3.0}). Thus if the solution of (\ref{3.0}) blows up
in finite time, then $\mathbb{E}u$ will do so.
The proof is complete. $\Box$

Theorem \ref{t3.2} shows that $\mathbb{E}u$ will be easier to blow up in finite time than the solution
of (\ref{3.0}), but do not give the blowup probability. Now, we study this interesting problem.
Let $u$ be a mild solution to (\ref{2.2}) in the sense of Walsh \cite{walsh1986} (given by the heat kernel). It follows Theorem \ref{t3.1} that the
solutions of (\ref{2.2}) will keep positive if the initial data are nonnegative. Furthermore, we want to know the probability of event that the solutions of (\ref{2.2})
blow up in finite time. It suffices to consider the equation (\ref{2.3}). For simplicity,
we only consider the case that $f(u)=u^p$. Following \cite{Hubook2018}, if $1<p\leq1+\frac{2}{d}$, then any nontrivial, nonnegative solution solutions of (\ref{2.3}) with $\sigma=0$
blows up in finite time. When $\sigma\neq0$, we have
\begin{theo}\lbl{t2.3} Assume that $1<p\leq 1+\frac{2}{pd}$.
The probability that the solution of (\ref{2.3}) blows up in
finite time is lower bounded by
$\Phi\left(\frac{\ln(\frac{1}{\varepsilon})-\frac{(p-1)\sigma^2 }{2}}{|\sigma|(p-1)}\right)$, where $0<\varepsilon\ll1$ is a fixed any small
constant;

Assume that $1+\frac{2}{pd}<p< 1+\frac{2}{d}$.
The probability that the solution of (\ref{2.3}) blows up in
finite time is lower bounded by
$\Phi\left(\frac{\ln(\frac{1}{\epsilon})-\frac{2(p-1)\sigma^2 T^\star}{2}}{|\sigma|(p-1)\sqrt{2T^\star}}\right)$, where $0<\epsilon\ll1$ is a fixed any small
constant and $T^\star$ satisfies (\ref{2.10}).
  \end{theo}

{\bf Proof.} It follows from Theorem \ref{t3.1} that $w(x,t)\geq0$ almost surely.
By using the properties of heat kernel,
we get
   \bess
w(x,t)&=&\int_{\mathbb{R}^d}K(t,x-y)u_0(y)dy
+\int_0^t\int_{\mathbb{R}^d}K(t-s,x-y)e^{-(p-1)\frac{\sigma^2}{2}s+(p-1)\sigma B_s}w^p(y,s)dyds\\
&=&:I_1(x,t)+I_2(x,t).
   \eess
We assume that the solution
remains finite for all finite $t$ almost surely and want to derive a contradiction.
We may assume without loss of generality that $u_0(x)\geq C_1>0$ for $|x|<1$ by the
assumption. A direct computation shows that
   \bess
I_1(x,t)&\geq&\frac{C_1}{(4\pi t)^{\frac{d}{2}}}\int_{B_1(0)}\exp\left(-\frac{|x|^2+|y|^2}{4t}\right)dy\nm\\
&\geq& \frac{C_1}{(4\pi t)^{\frac{d}{2}}}\exp\left(-\frac{|x|^2}{2t}\right)\int_{|y|\leq\frac{1}{\sqrt{t}}, \ y\in B_1(0)}
\exp\left(-\frac{|y|^2}{2t}\right)dy\nm\\
&\geq& \frac{C}{(4\pi t)^{\frac{d}{2}}}\exp\left(-\frac{|x|^2}{2t}\right)
  \lbl{2.4} \eess
for $t>1$ and $C>0$.

It is easy to see that
   \bess
I_2(x,t)&\geq& C_0\int_0^te^{-(p-1)\frac{\sigma^2}{2}s+(p-1)\sigma B_s}\left(\int_{\mathbb{R}^d}K(t-s,x-y)w(y,s)dy\right)^pds.
   \eess
Let
   \bess
G(t)=\int_{\mathbb{R}^d}K(t,x) w(x,t)dx.
   \eess
Then for $t>1$,
   \bes
G(t)&=&\int_{\mathbb{R}^d}I_1(x,t)K(t,x)dx+\int_{\mathbb{R}^d}I_2(x,t)K(t,x)dx\nm\\
&\geq&\frac{C_2}{t^{\frac{d}{2}}}+\int_0^te^{-(p-1)\frac{\sigma^2}{2}s+(p-1)\sigma B_s}\nm\\
&&\times\left(\int_{\mathbb{R}^d}\int_{\mathbb{R}^d}K(t,x)K(t-s,x-y)dxw(y,s)dy\right)^pds.
   \lbl{2.4}\ees
It follows from the estimates of \cite[pp 42]{Hubook2018} that
   \bess
\int_{\mathbb{R}^d}K(t,x)K(t-s,x-y)dx\geq C_3K(s,y)\frac{s^{\frac{d}{2}}}{t^{\frac{d}{2}}}.
   \eess
Hence, (\ref{2.4}) becomes
   \bess
G(t)\geq \frac{C_2}{t^{\frac{d}{2}}}+C_3\int_0^te^{-(p-1)\frac{\sigma^2}{2}s+(p-1)\sigma B_s}
\left(\frac{s^{\frac{d}{2}}}{t^{\frac{d}{2}}}\right)^pG^p(s)ds,
  \eess
where $C_3$ is a positive constant.
We can rewrite the above inequality as
   \bes
t^{pd/2} G(t)\geq C_2t^{d(p-1)/2}+C_3\int_0^ts^{pd/2}e^{-(p-1)\frac{\sigma^2}{2}s+(p-1)\sigma B_s}G^p(s)ds=:g(t).
\lbl{2.5}
   \ees
Then for $t>1$, we have
   \bess
&&g(t)\geq C_2t^{d(p-1)/2},\\
&&g'(t)\geq t^{pd/2}e^{-(p-1)\frac{\sigma^2}{2}t+(p-1)\sigma B_t}G^p(t)\geq t^{(p-p^2)d/2}e^{-(p-1)\frac{\sigma^2}{2}t+(p-1)\sigma B_t}g^p(t),
  \eess
which implies
   \bes
\frac{C_2^{1-p}}{p-1}t^{-\frac{d(p-1)^2}{2}}\geq \frac{1}{p-1}g^{1-p}(t)\geq C_3\int_t^Ts^{(p-p^2)d/2}e^{-(p-1)\frac{\sigma^2}{2}s+(p-1)\sigma B_s}ds\ \ {\rm for}\
T>t\geq1.
   \lbl{2.6}\ees

We first consider the case: $1<p\leq 1+\frac{2}{pd}$, i.e., $(p-p^2)d/2+1\geq0$.
Let $\{\mathcal {F}_t\}_{t\geq0}$ be the filtration generated by $\{B_t\}_{t\geq0}$. Due to
$\{e^{-\frac{\sigma^2}{2}t+\sigma B_t}\}_{t\geq0}$ is martingale with respect to $\{\mathcal {F}_t\}_{t\geq0}$, taking conditional expectation
and using Jensen's inequality, we have
   \bess
\frac{C_2^{1-p}}{p-1}t^{-\frac{d(p-1)^2}{2}}&\geq& C_3\int_t^Ts^{(p-p^2)d/2}
\mathbb{E}\left[e^{-(p-1)\frac{\sigma^2}{2}s+(p-1)\sigma B_s}|\mathcal {F}_t\right]ds\nm\\
&\geq&C_3\int_t^Ts^{(p-p^2)d/2}
\left(\mathbb{E}\left[e^{-\frac{\sigma^2}{2}s+\sigma B_s}|\mathcal {F}_t\right]\right)^{p-1}ds\nm\\
&\geq&C_3\int_t^Ts^{(p-p^2)d/2}ds
e^{-(p-1)\frac{\sigma^2}{2}t+(p-1)\sigma B_t}\nm.\eess
Observing that
  \bess
  \int_t^Ts^{(p-p^2)d/2}ds=\left\{\begin{array}{lll}
\frac{2}{(p-p^2)d+2}\left(T^{\frac{(p-p^2)d+2}{2}}-t^{\frac{(p-p^2)d+2}{2}}\right), &(p-p^2)d/2+1>0,\\[1.5mm]
   \ln T-\ln t, \ \ \ &(p-p^2)d/2+1=0,
    \end{array}\right.\eess
we gain
   \bes
&&\!\!\!\frac{C_2^{1-p}}{p-1}t^{-\frac{d(p-1)^2}{2}}\nm\\&\geq&\!\!\!
\left\{\begin{array}{lll}\!\!\!
\frac{2C_3}{(p-p^2)d+2}\left(T^{\frac{(p-p^2)d+2}{2}}-t^{\frac{(p-p^2)d+2}{2}}\right)
e^{-(p-1)\frac{\sigma^2}{2}t+(p-1)\sigma B_t}, &(p-p^2)d/2+1>0,\\[1.5mm]
\!\!\!C_3(\ln T-\ln t)e^{-(p-1)\frac{\sigma^2}{2}t+(p-1)\sigma B_t}, \ \ \ &(p-p^2)d/2+1=0.
\end{array}\right.
\lbl{2.7}\ees

Hence letting $T\to\infty$, we know that if probability of the inequality (\ref{2.7}) does not hold is equivalent to the probability of the
event that $\{\omega\in\Omega; \ e^{-(p-1)\frac{\sigma^2}{2}t+(p-1)\sigma B_t(\omega)}>\varepsilon\}$, where
$\varepsilon>0$ is any fixed. By using the fact that $B_t$ is a Guass process and for any fixed $t>0$,
$B_t$ obeys the Guass normal distribution $N(0,t)$, we have
       \bess
\mathbb{P}\left(e^{-(p-1)\frac{\sigma^2}{2}t+(p-1)\sigma B_t}>\varepsilon\right)=
\Phi\left(\frac{\ln(\frac{1}{\varepsilon})-\frac{(p-1)\sigma^2 t}{2}}{|\sigma|(p-1)\sqrt{t}}\right).
   \eess
It follows from the above discussion that we can take $t=1$ in above equality.

Now, we consider the case: $1+\frac{2}{pd}<p<1+\frac{2}{d}$. In this case, noting that $(p-p^2)d+2<0$ and $\frac{d(p-1)^2}{2}>-1+\frac{d(p-1)^2}{2}$,
we obtain that
 \bes
\frac{C_2^{1-p}}{p-1}t^{-\frac{d(p-1)^2}{2}}\geq \frac{2C_3}{(p^2-p)d-2}\left(t^{\frac{(p-p^2)d+2}{2}}-T^{\frac{(p-p^2)d+2}{2}}\right)
e^{-(p-1)\frac{\sigma^2}{2}t+(p-1)\sigma B_t}.
   \lbl{2.8}\ees
Letting $T\to\infty$, we get
 \bes
\frac{C_2^{1-p}}{p-1}t^{-\frac{d(p-1)^2}{2}}\geq \frac{2C_3}{(p^2-p)d-2}t^{\frac{(p-p^2)d+2}{2}}
e^{-(p-1)\frac{\sigma^2}{2}t+(p-1)\sigma B_t}.
   \lbl{2.9}\ees
For any fixed $\epsilon>0$, let $T^\star$ satisfy
   \bes
\frac{C_2^{1-p}}{p-1}(T^\star)^{-\frac{d(p-1)^2}{2}}< \frac{2C_3\epsilon}{(p^2-p)d-2}(T^\star)^{\frac{(p-p^2)d+2}{2}}.
   \lbl{2.10}\ees
Then the inequality (\ref{2.8}) does not hold with the probability $\mathbb{P}(A_{T^\star})$, where
   \bess
A_{T^\star}=\{\omega\in\Omega; \ e^{-(p-1)\frac{\sigma^2}{2}T^\star+(p-1)\sigma B_{T^\star}(\omega)}>\epsilon\}.
   \eess
Similar to the former case, we have
   \bess
\mathbb{P}(A_{T^\star})=\Phi\left(\frac{\ln(\frac{1}{\epsilon})-\frac{(p-1)\sigma^2 T^\star}{2}}{|\sigma|(p-1)\sqrt{T^\star}}\right).
  \eess

The proof is complete. $\Box$

\begin{remark}\lbl{r3.2}
The difference between the two cases $(p-p^2)d/2+1>0$ and $1+\frac{2}{pd}<p<1+\frac{2}{d}$ is that:
in the case $(p-p^2)d/2+1>0$, the constant $\varepsilon>0$ does not depend on the time; however, in
the case $1+\frac{2}{pd}<p<1+\frac{2}{d}$, the constant $\epsilon>0$ depends on the time and must
satisfy the inequality (\ref{2.10}). Consequently, in the case $(p-p^2)d/2+1>0$, the probability of
the event that  the solutions blow up in finite time is closed to $1$. But in the other case,
the probability has a certain distance with respect to $1$.

The linear multiplicative noise can be regarded as a perturbation and the profile of the
solution will keep together with the deterministic case. Thus we should be care of the
probability of the event that the solutions has the same properties as the deterministic case.
But if the noise is nonlinear multiplicative, the structure for the original equation will be changed, we can not deduce the same properties as the deterministic case in general.
\end{remark}

Lastly, we prove the probability of the event that blowup time of stochastic parabolic
equations large than or less than the deterministic case. Let $D\subset\mathbb{R}^d$. Consider the following stochastic parabolic equation
      \bes\left\{\begin{array}{lll}
   du=(\Delta u+G(u))dt+\kappa udB_t, \ \ \qquad &t>0,\ x\in  D,\\[1.5mm]
   u(x,0)=u_0(x), \ \ \ &\qquad\ \ x\in D,\\[1.5mm]
   u(x,t)=0, \ \ \ \ &t>0,\ x\in\partial D,
    \end{array}\right.\lbl{0.3}\ees
where $G:\mathbb{R}\to\mathbb{R}_+$ is locally Lipschitz and satisfies
    \bess
G(u)\geq C u^{1+\beta}\ \ \ {\rm for\ all}\ u>0,
  \eess
and $C,\beta,\kappa$ are positive numbers, $\{B_t\}_{t\geq0}$ is a standard
one-dimensional Brownian motion on a stochastic basis $(\Omega,\mathcal{F},\{\mathcal{F}_t\}_{t\geq0},\mathbb{P})$ and
$u_0:\ D\to\mathbb{R}_+$ is of class $C^2$ and not identically zero.
Dozzi and L\'{o}pez-Mimbela
\cite{DL2010} obtained the probability that the solution of (\ref{0.3})
blows up in finite time is lower bounded by
$\int_{\frac{1}{\beta}u(\phi,0)^{-\beta}}^{+\infty}h(y)dy$, where
   \bess
h(y)=
\frac{(\kappa^2\beta^2y/2)^{(2\lambda_1+\kappa^2)/\kappa^2\beta}}
{y\Gamma((2\lambda_1+\kappa^2)/(\kappa^2\beta))}\exp\left(-\frac{2}{\kappa^2\beta^2y}\right),\ \
u(\phi,0)=\int_Du_0(x)\phi(x)dx,
    \eess
where $\lambda_1>0$ is the first eigenvalue of the Laplacian on $D$, and
$\phi$ is the corresponding eigenfunction normalized so that $\|\phi\|_{L^1}=1$.
It is not hard to prove that if let $v(x,t)=e^{-\kappa B_t}u(x,t)$, then $v$ satisfies
      \bes\left\{\begin{array}{lll}
\frac{\partial}{\partial t}v(x,t)=\Delta v(x,t)-\frac{\kappa^2}{2}v(x,t)+e^{-\kappa B_t}
G(e^{\kappa B_t}v(x,t)), \ \ \qquad &t>0,\ x\in  D,\\[1.5mm]
   v(x,0)=u_0(x), \ \ \ &\qquad\ \ x\in D,\\[1.5mm]
   v(x,t)=0, \ \ \ \ &t>0,\ x\in\partial D.
    \end{array}\right.\lbl{0.4}\ees
%By using Jensen's inequality, one can prove that
 %  \bess
%\frac{d}{dt}v(\phi,t)\geq-\left(\lambda_1+\frac{\kappa^2}{2}\right)v(\phi,t)
%+e^{\kappa\beta B_t}v(\phi,t)^{1+\beta}.
 %  \eess
Similar to the proof of Theorem \ref{t3.1}, we can prove the solutions of
(\ref{0.4}) will keep positive. Following the method of Theorem \ref{t2.3}, one
can give a different probability from \cite{DL2010} of the event that the solutions blow up in finite time.
In paper \cite{DL2010}, the blowup time is obtained, that is,
   \bess
\tau:=\inf\left\{t\geq0 \ \Big| \ \int_0^te^{-(\lambda_1+\kappa^2/2)\beta s+\kappa\beta B_s}ds
\geq\frac{1}{\beta}u(\phi,0)\right\}.
   \eess
It is easy to see that when $\kappa=0$, $\tau$ becomes the blowup time of
deterministic case.  Assume $T^*$ satisfies
   \bess
\int_0^{T^*}e^{-\lambda_1\beta s}ds
=\frac{1}{\beta}u(\phi,0).
    \eess
Now we want to prove $\mathbb{P}(\tau>T^*)$. It follows from the definition of
$T^*$, we have
   \bess
&&\mathbb{P}(\tau>T^*)\\
&=&\mathbb{P}\left(\tau>T^*,\,\int_0^{T^*}e^{-\lambda_1\beta s}
(1-e^{-\kappa^2\beta s/2+\kappa\beta B_s})ds=
\int_{T^*}^\tau e^{-(\lambda_1+\kappa^2/2)\beta s+\kappa\beta B_s}ds\right)\\
&=&\mathbb{P}\left(\int_0^{T^*}e^{-\lambda_1\beta s}
(1-e^{-\kappa^2\beta s/2+\kappa\beta B_s})ds>0\right)\\
&=&\mathbb{P}\left(\int_0^{T^*}e^{-(\lambda_1+\kappa^2/2)\beta s+\kappa\beta B_s}ds<
\frac{1}{\lambda_1\beta}(1-e^{-\lambda_1\beta T^*})\right).
   \eess
It follows from the results of \cite{Yor1992} that the random variable
$\int_0^{T^*}e^{-(\lambda_1+\kappa^2/2)\beta s+\kappa\beta B_s}ds$ has a
probability law and we denote by $P^{T^*}$. Then we have
   \bess
\mathbb{P}(\tau>T^*)=P^{T^*}(\frac{1}{\lambda_1\beta}(1-e^{-\lambda_1\beta T^*})).
  \eess
Combining the above discussion, we have the following result.
\begin{theo}\lbl{t2.4} The probability that the blowup time of (\ref{0.3}) is
larger than the deterministic case (i.e., (\ref{0.3}) with $\kappa=0$) is
$P^{T^*}(\frac{1}{\lambda_1\beta}(1-e^{-\lambda_1\beta T^*}))$.
  \end{theo}

\medskip

\noindent {\bf Acknowledgment} The first author was supported in part
by NSFC of China grants 11771123. The authors thanks Prof. Feng-yu Wang for discussing
this manuscript.

 \end{document}